\documentclass[12pt]{article}
\usepackage{amsfonts,amsmath,amsxtra}
\usepackage{amssymb, amscd}
\def\hybrid{\topmargin 0pt      \oddsidemargin 0pt
        \headheight 0pt \headsep 0pt
        \textwidth 16.5cm
        \textheight 23cm
        \marginparwidth 0.0in
        \parskip 5pt plus 1pt   \jot = 1.5ex}
\catcode`\@=11
\def\marginnote#1{}
\newcount\hour
\newcount\minute
\newtoks\amorpm
\hour=\time\divide\hour by60
\minute=\time{\multiply\hour by60 \global\advance\minute by-\hour}
\edef\standardtime{{\ifnum\hour<12 \global\amorpm={am}%
        \else\global\amorpm={pm}\advance\hour by-12 \fi
        \ifnum\hour=0 \hour=12 \fi
      \number\hour:\ifnum\minute<10 0\fi\number\minute\the\amorpm}}
\edef\militarytime{\number\hour:\ifnum\minute<10 0\fi\number\minute}

\def\draftlabel#1{{\@bsphack\if@filesw {\let\thepage\relax
   \xdef\@gtempa{\write\@auxout{\string
      \newlabel{#1}{{\@currentlabel}{\thepage}}}}}\@gtempa
   \if@nobreak \ifvmode\nobreak\fi\fi\fi\@esphack}
        \gdef\@eqnlabel{#1}}
\def\@eqnlabel{}
\def\@vacuum{}
\def\draftmarginnote#1{\marginpar{\raggedright\scriptsize\tt#1}}

\def\draft{\oddsidemargin -0.1truein
        \def\@oddfoot{\sl preliminary draft \hfil
        \rm\thepage\hfil\sl\today\quad\militarytime}
        \let\@evenfoot\@oddfoot \overfullrule 3pt
        \let\label=\draftlabel
        \let\marginnote=\draftmarginnote
\def\@eqnnum{{\rm (\theequation)}
\rlap{\kern\marginparsep\tt\@eqnlabel}%
\global\let\@eqnlabel\@vacuum}  }
\newcommand{\RR}{{\mathbb{R}}}
\newcommand{\CC}{{\mathbb{C}}}
\newcommand{\ZZ}{{\mathbb{Z}}}

\newfont{\Bbbb}{msbm7 scaled 1\@ptsize00}
\newcommand{\zs}{\raise-1pt\hbox{$\mbox{\Bbbb Z}$}}

\font\sevenmsa=msam6 %scaled 1\@ptsize00
\newfam\msafam
\textfont\msafam=\sevenmsa
\def\hexnumber@#1{\ifnum#1<10 \number#1\else
\ifnum#1=10 A\else\ifnum#1=11 B\else\ifnum#1=12 C\else
\ifnum#1=13 D\else\ifnum#1=14 E\else\ifnum#1=15 F\fi\fi\fi\fi\fi\fi\fi}
\def\msa@{\hexnumber@\msafam}
\def\llcorner{\delimiter"4\msa@78\msa@78 }
\def\lrcorner{\delimiter"5\msa@79\msa@79 }
\mathchardef\blacktriangleright="3\msa@49
\mathchardef\blacktriangleleft="3\msa@4A
\font\tenmsb=msbm10 scaled 1\@ptsize00
\newfam\msbfam
\textfont\msbfam=\tenmsb
\scriptfont\msbfam=\tenmsb
\newdimen\linethick  \linethick=0.4pt
\newdimen\hboxitspace    \hboxitspace=5pt
\newdimen\vboxitspace    \vboxitspace=5pt

\def\fr#1{%
\be\new
\vcenter{
\hrule height\linethick
           \hbox{\vrule width\linethick
                 \kern\hboxitspace
                 \vbox{\kern\vboxitspace
                       \hbox{$\begin{array}{c}\displaystyle#1
          \end{array}$}%
                       \kern\vboxitspace}%
                 \kern\hboxitspace
                 \vrule width\linethick}%
           \hrule height\linethick}%
\ee}
\newdimen\Squaresize \Squaresize=14pt
\newdimen\Thickness \Thickness=0.5pt

\def\Square#1{\hbox{\vrule width \Thickness
   \vbox to \Squaresize{\hrule height \Thickness\vss
      \hbox to \Squaresize{\hss#1\hss}
   \vss\hrule height\Thickness}
\unskip\vrule width \Thickness}
\kern-\Thickness}

\def\Vsquare#1{\vbox{\Square{$#1$}}\kern-\Thickness}

\def\numberbysection{\@addtoreset{equation}{section}
        \def\theequation{\thesection.\arabic{equation}}}
\numberbysection

\renewcommand{\theequation}{\thesection.\arabic{equation}}
\def\titlepage{\@restonecolfalse\if@twocolumn\@restonecoltrue\onecolumn
     \else \newpage \fi \thispagestyle{empty}\c@page\z@
        \def\thefootnote{\fnsymbol{footnote}} }

\def\endtitlepage{\if@restonecol\twocolumn \else  \fi
        \def\thefootnote{\arabic{footnote}}
        \setcounter{footnote}{0}}  %\c@footnote\z@ }
\relax

\hybrid
\makeatletter
\newdimen\normalarrayskip            % skip between lines
\newdimen\minarrayskip               % minimal skip between lines
\normalarrayskip\baselineskip
\minarrayskip\jot
\newif\ifold             \oldtrue            \def\new{\oldfalse}
\def\arraymode{\ifold\relax\else\displaystyle\fi}%mode of array enrties
\def\eqnumphantom{\phantom{(\theequation)}} % ight phantom in eqnarray
\def\@arrayskip{\ifold\baselineskip\z@\lineskip\z@
     \else
     \baselineskip\minarrayskip\lineskip1\baselineskip\fi}

\def\@arrayclassz{\ifcase \@lastchclass \@acolampacol \or
\@ampacol \or \or \or \@addamp \or
   \@acolampacol \or \@firstampfalse \@acol \fi
\edef\@preamble{\@preamble
  \ifcase \@chnum
     \hfil$\relax\arraymode\@sharp$\hfil
     \or $\relax\arraymode\@sharp$\hfil
     \or \hfil$\relax\arraymode\@sharp$\fi}}

\def\@array[#1]#2{\setbox\@arstrutbox=\hbox{\vrule
     height\arraystretch \ht\strutbox
     depth\arraystretch \dp\strutbox
width\z@}\@mkpream{#2}\edef\@preamble{\halign \noexpand\@halignto
\bgroup \tabskip\z@ \@arstrut \@preamble \tabskip\z@ \cr}%
\let\@startpbox\@@startpbox \let\@endpbox\@@endpbox
  \if #1t\vtop \else \if#1b\vbox \else \vcenter \fi\fi
  \bgroup \let\par\relax
  \let\@sharp##\let\protect\relax
  \@arrayskip\@preamble}
\def\eqnarray{\stepcounter{equation}%
              \let\@currentlabel=\theequation
              \global\@eqnswtrue
              \global\@eqcnt\z@
              \tabskip\@centering              %formulae  centering
              \let\\=\@eqncr
              $$%
            \halign to \displaywidth  \bgroup
             \eqnumphantom \@eqnsel
      \hskip\@centering                               %right tab%
    $\displaystyle  \tabskip\z@ {##}$%
    &\global\@eqcnt\@ne \hskip 2\arraycolsep
         $ \displaystyle  \arraymode{##}$\hfil
    &\global\@eqcnt\tw@ \hskip 2\arraycolsep
         $\displaystyle\tabskip\z@{##}$\hfil
         \tabskip\@centering
    &{##}\tabskip\z@\cr}
\makeatother
\newtheorem{te}{Theorem}[section]%Usage:\begin{te}Statement\end{te}
\newtheorem{de}{Definition}[section]
\newtheorem{prop}{Proposition}[section]           %  ETC ...

\newtheorem{lem}{Lemma}[section]

\newcommand{\beq}[1]{\begin{equation}\label{#1}}
\newcommand\eeq{\end{equation}}
\newcommand\bqa{\begin{eqnarray}}
\newcommand\eqa{\end{eqnarray}}

\def\be{\begin{eqnarray}\new\begin{array}{cc}}
\def\ee{\end{array}\end{eqnarray}}

\def\beq{\begin{equation}}
\def\eeq{\end{equation}}
\def\bse{\begin{subequations}}                %%%SUBEQUATIONS
\def\ese{\end{subequations}}
\def\bp{\begin{pmatrix}}
\def\ep{\end{pmatrix}}

\def\h{\hbar}
\def\i{\imath}
\renewcommand{\theequation}{\thesection.\arabic{equation}}

\def\square{\hfill{\vrule height6pt width6pt            %Black
depth1pt} \break \vspace{.01cm}}                        %square
\def\d{\partial}
\def\stack#1#2{\raise0.7pt\hbox{$\mathrel{\mathop{#2}\limits^{#1}}$}}
\def\tr{\triangleright}
\def\tl{\triangleleft}
\def\sem{\mathsurround=0pt \raise1pt
\hbox{$\scriptscriptstyle>\!\!$}\:\!\!\tl}
\def\mes{\mathsurround=0pt \tr\!\:\!\raise0.8pt
\hbox{$\scriptscriptstyle\!\!<$}\,}
\def\]{\mathsurround=0pt ]\raise-2pt\hbox{$_\ast$}}
\def\al{\alpha}

\def\<{\langle}
\def\>{\rangle}
\def\ov{\overline}
\def\wt{\widetilde}
\def\wh{\widehat}
\def\vk{\varkappa}
\def\frak{\mathfrak}

\def\N{{\scriptscriptstyle N}}
\def\ts#1#2{{\textstyle\frac{#1}{#2}}}
\def\g{\frak g}

\def\we{\raise-1pt\hbox{$\,\stackrel{\wedge}{,}\,$}}
\def\t{{\rm tr}\,}
\def\pr {\partial}
\def\ay{A}
\def\by{B}
\def\cy{C}
\def\dy{D}

0
\begin{document}
\def\t{\theta}
\def\T{\Theta}
\def\w{\omega}
\def\ov{\overline}
\def\a{\alpha}
\def\b{\beta}
\def\g{\gamma}
\def\s{\sigma}
\def\l{\lambda}
\def\wt{\widetilde}

%\draft
%\begin{flushright}
%monopole.tex
%\hfill{\normalsize ITEP-TH-1/04}\\ [10mm]\end{flushright}

%\draft                             %SWITCH ON/OFF DRAFT VERSION%
\thispagestyle{empty}
\begin{center}

\phantom.
\bigskip%\bigskip\bigskip\bigskip\bigskip\bigskip
{\hfill{\normalsize ITEP-TH-30/04}\\
\hfill{\normalsize TCD-MATH-04-15}\\
\hfill{\normalsize HMI-04-04}\\
%\hfill{\normalsize MPIM-....}\\
[15mm]\Large\bf
 On a class of representations of the  Yangian  \\  and
  moduli space of monopoles}

\vspace{1cm}

\bigskip\bigskip
{\large A. Gerasimov}
\\ \bigskip
{\it
Institute for Theoretical \& Experimental Physics, 117259, Moscow,
Russia}\\ {\it
 Department of Pure and Applied Mathematics, Trinity
College, Dublin 2, Ireland } \\ {\it Hamilton
Mathematics Institute, TCD, Dublin 2, Ireland}\\
\bigskip
{\large S. Kharchev\footnote{E-mail: kharchev@itep.ru}},\\
\bigskip
{\it Institute for Theoretical \& Experimental Physics, 117259, Moscow,
Russia}\\
\bigskip
{\large D. Lebedev\footnote{E-mail: lebedev@mpim-bonn.mpg.de}}
\\ \bigskip
{\it Institute for Theoretical \& Experimental Physics, 117259, Moscow,
Russia} and {\it Max-Planck-Institut fur Mathematik, Vivatsgasse 7, D-53111
Bonn, Germany},\\
\bigskip
{\large S. Oblezin}
\footnote{E-mail: Sergey.Oblezin@itep.ru}\\ \bigskip
{\it
Institute for Theoretical \& Experimental Physics, 117259, Moscow,
Russia}\\
\end{center}

%\vspace{1cm}

%\bigskip \bigskip
%{\large A. Gerasimov}\footnote{
%{\it
%Institute for Theoretical \& Experimental Physics, 117259, Moscow,
%Russia}}\footnote{{\it
% Department of Pure and Applied Mathematics, Trinity
%College, Dublin 2, Ireland }}\footnote{{\it Hamilton
%Mathematics Institute, TCD, Dublin 2, Ireland}},
%\bigskip\bigskip
%{\large S. Kharchev{\small $^1$}, D.
%Lebedev}{\small $^1$}\footnote{{\it Max-Plank-Institut f\"{u}r Mathematik, Vivatsgasse 7, D-53111
% Bonn, Germany}},
%{\large and S. Oblezin{\small
%$^1$}}
%\end{center}

\vspace{0.5cm}

\begin{abstract}
\noindent

A new class of infinite dimensional  representations of the Yangians
$Y(\frak{g})$ and $Y(\frak{b})$ corresponding to a complex semisimple algebra
$\frak{g}$ and its Borel subalgebra $\frak{b}\subset\frak{g}$ is constructed.
It is based on the generalization of the Drinfeld realization of
$Y(\frak{g})$, $\frak{g}=\frak{gl}(N)$  in terms of quantum minors
to the case of an arbitrary  semisimple Lie algebra $\frak{g}$.
The Poisson geometry associated with the constructed representations
is described. In particular  it is shown that the underlying
symplectic leaves  are isomorphic to the moduli
spaces of $G$-monopoles  defined as the components of the space of based
maps of $\mathbb{P}^1$ into the generalized flag manifold $X=G/B$. Thus the
constructed representations of the Yangian may be considered as a
quantization of the moduli space of the monopoles.
\end{abstract}

%\end{titlepage}
\clearpage \newpage
\setcounter{footnote}0

%\draft                             %SWITCH ON/OFF DRAFT VERSION%

%\tableofcontents
\normalsize
\section{Introduction}
The Yangian $Y(\frak{g})$ for a semisimple complex Lie algebra
$\frak{g}$ was introduced by Drinfeld as a certain deformation of
$U(\frak{g}[t])$ as a Hopf algebra \cite{Dr1}, \cite{Dr2},
\cite{Dr3} (see for  recent review \cite{MNO}, and \cite{M}).
Recently   a construction of the special class of
infinite-dimensional representations of $Y(\frak{gl}(N))$ and
$Y(\frak{sl}(N))$ based on the generalization of the
Gelfand-Zetlin construction was introduced  in \cite{GKL1} (see also
\cite{GKL2}).  In this paper we generalize  this construction to  $Y(\frak{g})$ for
an arbitrary semisimple Lie algebra $\frak{g}$. This
generalization is based on  the proposed generalization of the
Drinfeld realization of $Y(\frak{gl}(N))$ in terms of quantum
minors to the case of $Y(\frak{g})$ for an arbitrary semisimple
Lie algebra $\frak{g}$. We also construct  representations of
$Y(\frak{b})$ where $\frak{b}\subset\frak{g}$ is a Borel
subalgebra of $\frak{g}$. One should note that not all of the considered
representations of $Y(\frak{b})$ may be lifted to the
representations of $Y(\frak{g})$. We also  describe  the
Poisson geometry behind the constructed representation by defining
explicitly the symplectic leaves of the  classical versions of the
Yangians $Y(\frak{g})$ and $Y(\frak{b})$. The proposed description of
the symplectic leaves  reveals a deep connection with the moduli
space of $G$-monopoles, $Lie(G)=\frak{g}$. The symplectic leaves
as symplectic manifolds turn out to be open parts of the moduli
space  of monopoles  with maximal symmetry breaking supplied with
the symplectic structure introduced in \cite{AH} for $G=SU(2)$, in
\cite{B2} for  $G=SU(N)$, and  in \cite{FKMM}  for an arbitrary
semisimple Lie group $G$. Thus our construction of the Yangian
representations
 can be considered as a quantization of the  moduli spaces of
$G$-monopoles.

The results of this paper support the strong connection between
quantum integrable systems and problems of the quantization of
various moduli spaces. The demonstrated connection between the Yangian and the
quantization of the moduli space of the monopoles is a particular
example of this deep relation. We plan to consider its  implications
  to the theory of quantum integrable theories elsewhere. Let us
remark that the connection between the Atiyah-Hitchin symplectic structure on
the moduli space of the $SU(2)$ monopoles \cite{AH} with some  particular
integrable systems was noted previously in %\cite{JJJ},
\cite{Van1}, \cite{Van2}.

Finally note that the explicit construction of the representations
of $Y(\frak{g})$ discussed below appears to be similar
to the constructions of the representations of a class of  elliptic algebras
proposed in \cite{FO1}, \cite{FO2}. Nevertheless, the main result
(Theorem \ref{MainTH}) seems to be new.

The plan of the paper is as follows. In Section 2 we provide  various
descriptions of the $Y(\frak{g})$ in terms of the generators and
relations. The construction  of $Y(\frak{g})$  in terms of the generators
$A_i(u),B_i(u),C_i(u)$ for a general Lie algebra is proposed.
In Section 3, we describe a particular class of representations  of
the $Y(\frak{g})$ and $Y(\frak{b})$ and give an explicit realization of
its generators in terms of  difference operators.  The main
result is formulated in the Theorem \ref{MainTH}.
In Section 4 we discuss the  underlying Poisson geometry and provide a
description of the corresponding  symplectic leaves of the classical
counterpart of the  Yangian. It appears that there is an isomorphism
between the open part of the symplectic leaves for $Y(\frak{b})$ and the
open part of  the moduli space of the $G$-monopoles with the maximal
symmetry breaking.

{\em Acknowledgments}: The  authors are grateful to A. Levin and A. Rosly
for useful discussions and to M. Finkelberg for the explanation of the results
of \cite{FKMM}, \cite{FM}. The research was partly supported by grants
CRDF RM1-2545; INTAS 03-513350; grant NSh 1999.2003.2 for support of
scientific schools, and by grants RFBR-040100646 (A. Gerasimov, D. Lebedev),
and RFBR-040100642 (S. Kharchev, S. Oblezin). The research of
A. Gerasimov was also partly supported by SFI Basic Research Grant.
D. Lebedev would like to thank the Max-Planck-Institut f\"{u}r Mathematik
for financial support and hospitality. S. Oblezin is deeply grateful to the
Independent University of Moscow for the support.

\section{The various presentations of $Y(\frak{g})$}
We start with  the definition of the Yangian for a semisimple
Lie algebra $\frak{g}$ due to Drinfeld  \cite{Dr1}, \cite{Dr2}
(see also \cite{CP}) in the form given in  \cite{KT}.

Let
 $\mathfrak{h}\subset\mathfrak{b}\subset\mathfrak{g}$ be
a simple finite dimensional Lie algebra $\frak{g}$ of rank $\ell$
over $\CC$ with a Cartan subalgebra $\frak{h}$ and a
 Borel subalgebra $\frak{b}$. Let $a=||a_{ij}||,
 i,j=1,\ldots,\ell$ be the Cartan matrix of $\frak g$,
$\Gamma$ be the set of vertices of the Dynkin diagram of
$\mathfrak{g}$, $\{\alpha_i \in\frak{h}^{*}, i\in \Gamma\}$ be
the set of simple roots and $\{\alpha_{i}^{\vee}, i\in\Gamma\}$ be the
set of the corresponding co-roots  ($a_{ij}=\alpha^{\vee}_{i}(\alpha_{j})$).
There exist positive integers $d_1,\ldots,d_\ell$ such that the
matrix $||d_{i}a_{ij}||$ is symmetric. Define the invariant bilinear form
on $\frak{h}^{*}$ by $(\alpha_i ,\alpha_j)=d_i a_{ij}$, then
$a_{ij}=\frac{2(\alpha_i,\alpha_j)}{(\alpha_i,\alpha_i)}$.

It is convenient to define the generators of the Yangian in terms
of  the generating  series $H_i (u)\,,
E_i  (u),\,$ and $ F_i(u)\,,\;i\in\Gamma$:
\be\label{int7}
H_i(u)=1+\sum_{s=0}^\infty H_i^{(s)}u^{-s-1}\;,\\
E_i(u)=\sum_{s=0}^\infty E_i^{(s)}u^{-s-1}\,,\ \ \ \ \ \
F_i(u)=\sum_{s=0}^\infty F_i^{(s)}u^{-s-1}\;.
\ee

\begin{de} The Yangian $Y(\frak{g})$ is the associative algebra with
 generators $H^{(s)}_i, E^{(s)}_i, F^{(s)}_i$, $i\in \Gamma; s=0,1,\ldots$
and the following defining relations
\be\label{cr1}
[H_i(u),H_j(v)]=0\;,
\ee
\be\label{cr2}
[H_i(u),E_j(v)]=-\frac{\imath\hbar}{2}\, (\alpha_i,\alpha_j)\,
\frac{[H_i(u),E_j(u)-E_j(v)]_{+}}{u-v}\,,\\

[H_i(u),F_j(v)]=\frac{\imath\hbar}{2}\, (\alpha_i,\alpha_j)\,
\frac{[H_i(u),F_j(u)-F_j(v)]_{+}}{u-v},
\ee
\be\label{cr3}
[E_i(u),F_j(v)]=-\imath\hbar\,\frac{H_i(u)-H_i(v)}
{u-v}\,\delta_{i,j}\,,
\ee

\be\label{cr4}
\hspace{-0.5cm}
[E_i(u),E_i(v)]=
-\frac{\imath\hbar}{2}(\alpha_i,\alpha_i)\frac{(E_i(u)-E_i(v))^2}{u-v}\,,\\
\hspace{-0.5cm}
[F_i(u),F_i(v)]=
\frac{\imath\hbar}{2}(\alpha_i,\alpha_i)\frac{(E_i(u)-E_i(v))^2}{u-v}\,,\\
\hspace{-0.5cm}
[E_i(u),E_j(v)]=-\frac{\imath\hbar}{2}
(\alpha_i,\alpha_j)\frac{[E_i(u),E_j(u)-E_j(v)]_+}{u-v}-
\frac{[E_i^{(0)},E_j(u)-E_j(v)]}{u-v}\,,\\
\hspace{-0.5cm}
[F_i(u),F_j(v)]=\frac{\imath\hbar}{2}
(\alpha_i,\alpha_j)\frac{[F_i(u),F_j(u)-F_j(v)]_+}{u-v}-
\frac{[F_i^{(0)},F_j(u)-F_j(v)]}{u-v}\,,\\
i\neq j,\,\, a_{ij}\neq 0;
\ee
%\newpage
\be\label{cr5}
\sum_{\sigma\in\frak S_n} [E_i(u_{\sigma(1)}),[E_i(u_{\sigma(2)}),
\ldots,[E_i(u_{\sigma(n)}),E_j(v)]\ldots]]=0\;,\\
\sum_{\sigma\in\frak S_n}[F_i(u_{\sigma(1)}),[F_i(u_{\sigma(2)}),
\ldots,[F_i(u_{\sigma(n)}),F_j(v)]\ldots]]=0\;,
\\ i\neq j,\; n=1-a_{ij}\,,
\ee
where  $[a,b]_+:= ab+ba$.
\end{de}
Denote $Y(\frak{b})\subset Y(\frak{g})$ the subalgebra generated
by $H_i(u), E_i(u),\,i\in\Gamma$.

There is another closely related set of  generators of the Yangian
which were given
by Drinfeld \cite{Dr3} in the case $Y(\frak{sl}(\ell+1))$.
Below we propose the generalization of this description to
$Y(\frak{g})$ for an
arbitrary semisimple Lie algebra $\frak{g}$. In sequel we will use the
following convention: $\prod\limits_{s=j}^k f_s=1$, for any $f_s$ if $k\le j$
which help to write down the formulas in a compact way.
\begin{lem}
Any generation series $H_i(u)$ of the form (\ref{int7}) can be
represented in the following form
\be\label{rf3}
\hspace{2cm}
H_i(u)=\frac{\prod\limits_{s\neq i}\prod\limits_{r=1}^{-a_{si}}
A_s\big(u-\ts{\i\h}{2}(\al_i+r\al_s,\al_s)\big)}
{A_i(u)A_i(u-\ts{\i\h}{2}(\al_i,\al_i))}\,,\ \ \  i=1,\ldots,\ell.
\ee
where $A_i(u)$, $i=1,\ldots,\ell$ are  formal series
$A_i
(u)=1+\sum\limits_{s=0}^{\infty}A_i^{(s)}u^{-s-1}$.
\end{lem}
Proof: Expanding the left and right hand sides
of (\ref{rf3}) in powers $u^{-1}$ and equating the coefficients in
$u^{-r-1}$ one can check that
\be\label{rf3'}
H_i^{(r)}=\sum\limits_{s=1}^\ell a_{si} A_s^{(r)}+f(A_i^{(0)},\cdots,
A_i^{(r-1)}).
\ee
Since the Cartan matrix is invertible, the statement follows by the
induction procedure \square

%The examples  of the explicit expressions $A_i(u)$
%through $H_i(u)$ can be found in Appendix A.

Let us introduce  the generating series $C_i(u)$ and $B_i(u)$:
\be\label{rf8}
B_i(u)=d_i^{1/2}A_i(u)E_i(u),\;\;\;C_i(u)=
d_i^{1/2}F_i(u)A_i(u)\,, \ \ \ i\in\Gamma\,,
\ee
where $d_i=(\al_i,\al_i)/2$. Due to (\ref{rf3}), the change of the generators
$H_i(u),E_i(u),F_i(u)$ to $A_i(u),C_i(u),B_i(u)$ is invertible.
It is  convenient to  introduce an additional set of  series $D_i(u)$.
For any $i\in\Gamma$ let
\be\label{rf12}
D_i(u)=A_i(u)H_i(u)+C_i(u)A_i^{-1}(u)B_i(u)\,,
\ee
where $H_i(u)$ are expressed in terms of $A_i(u)$  in (\ref{rf3}). Straightforward
calculations lead to the following statement.
\begin{prop}
For any $i=1,\ldots,\ell$ the series
$A_i(u),B_i(u),C_i(u),D_i(u)$ satisfy the relations
\be\label{y6}
[A_i(u),A_j(v)]=0\,,
\ee
\be
\hspace{2cm}
[A_i(u),B_j(v)]=[A_i(u),C_j(v)]=0\,,\hspace{1.2cm}(i\neq j)\,,\\
\hspace{2cm}
[B_i(u),B_j(v)]=[C_i(u),C_j(v)]=0\,,\hspace{1cm}(a_{ij}=0,\,i\neq j)\,,
\hspace{-1cm}\\

[B_i(u),B_i(v)]=[C_i(u),C_i(v)]=0\,,\\
\hspace{3cm}
[B_i(u),C_j(v)]=0\,,\hspace{2cm}(i\neq j)\,,
\ee
\be\label{y7}
(u-v)[\ay_i(u),\by_i(v)]=\i\h d_i
\big(\by_i(u)\ay_i(v)-\by_i(v)\ay_i(u)\big)\,,
\ee
\be\label{y9}
(u-v)[\ay_i(u),\cy_i(v)]=\i\h d_i
\big(\ay_i(u)\cy_i(v)-\ay_i(v)\cy_i(u)\big)\,,
\ee
\be\label{y11}
(u-v)[B_i(u),C_i(v)]\,=\,\i\h d_i\big(A_i(u)D_i(v)-A_i(v)D_i(u)\big)\,,
\ee
\be\label{y8}
(u-v)[\by_i(u),\dy_i(v)]=\i\h d_i
\big(\by_i(u)\dy_i(v)-\by_i(v)\dy_i(u)\big)\,,
\ee
\be\label{y10}
(u-v)[\cy_i(u),\dy_i(v)]=\i\h d_i
\big(\dy_i(u)\cy_i(v)-\dy_i(v)\cy_i(u)\big)\,,
\ee
\be\label{y12}
(u-v)[A_i(u),D_i(v)]=\i\h d_i\big(B_i(u)C_i(v)-B_i(v)C_i(u)\big)\,.
\ee
\end{prop}
In addition, the following analog of the quantum determinant relation holds:
\be\label{rm13}
\ay_i(u)\dy_i(u-\ts{\i\h}{2}(\al_i,\al_i))-\by_i(u)\cy_i
(u-\ts{\i\h}{2}(\al_i,\al_i))\,=\\
%\dy_i(u-\ts{\i\h}{2}(\al_i,\al_i))\ay_i(u)-
%\by_i(u-\ts{\i\h}{2}(\al_i,\al_i))\cy_i(u)\,=\\
%\dy_i(u)\ay_i(u-\ts{\i\h}{2}(\al_i,\al_i))-
%\cy_i(u)\by_i(u-\ts{\i\h}{2}(\al_i,\al_i))=\\
%\ay_i(u-\ts{\i\h}{2}(\al_i,\al_i))\dy_i(u)-
%\cy_i(u-\ts{\i\h}{2}(\al_i,\al_i))\by_i(u)=\\
\prod\limits_{s\neq i}\prod\limits_{r=1}^{-a_{si}}
\ay_s(u-\ts{\i\h}{2}(\al_i+r\al_s,\al_s))\,.
\ee
The described set of relations between $A_i(u),B_i(u),C_i(u)$ and
$D_i(u)$ is a  generalization of
the $Y(\frak{sl}(\ell+1))$ relations of \cite{NT1} to the case
of $Y(\frak{g})$ with an arbitrary simple Lie algebra $\frak{g}$.

By (\ref{y6}), the coefficients of the generating series
$A_i(u),\,i=1,\ldots,\ell$ define a commutative subalgebra of  $Y(\frak{g})$.
This property will be important in the construction of the representations
of the Yangian in the next section.

\section{Construction of the representations of $Y(\frak{g})$ and
$Y(\frak{b})$}

In this section the explicit construction of the representations of
the Yangian in terms of difference operators is given.
The explicit description of the representation of the Yangian in terms of
difference/differential operators acting in some functional space  is
based  on the choice of a commutative subalgebra. The elements of this
subalgebra  act in this representation by
multiplication on  functions. The proposed
construction uses the subalgebra generated by
$A_i(u)$ as a distinguished commutative subalgebra. However, we start with
explicit description of the resulting representation and then we make
some comments how it could be derived starting with the representation
of commutative subalgebra generated by $A_i(u)$ and using
 the commutation relations between $A_i(u),B_i(u),$ and $C_i(u)$ described in
Section 2.

Let us introduce a set of variables
$\{\gamma_{i,k}\,\,;i\in\Gamma; \,k=1,\ldots,m_i\}$,
 where $m_i$ are  arbitrary positive  integer numbers
 and let $\cal M$ be the space of meromorphic functions in these variables.
Let us define the following difference operators acting on $\cal
M$:\bqa\label{shift}
\beta_{i,k}=e^{\frac{\i\hbar}{2}(\al_i,\al_i)\d_{\gamma_{i,k}}}.
\eqa
where $\d_{\gamma_{i,k}}:=\frac{\pr}{\pr \gamma_{i,k}}$ is a
differentiation over $\gamma_{i,k}$.
It is useful to arrange the variables into the set of  polynomials
 of the   formal variable $u$  of degrees $m_i,\,i\in\Gamma$
\be\label{mn1}
\hspace{1cm}
P_i(u)=\prod\limits_{p=1}^{m_i}(u-\gamma_{i,p})\,,\ \ \ \ i\in\Gamma.
\ee
%We shall use the notation
%\be\label{mn2}
%P_{i}'(\gamma_{i,k}):=\frac{\pr P_i(u)}{\pr u}\Big|_{u=\gamma_{i,k}}=
%\prod\limits_{p\neq k}(\gamma_{i,k}-\gamma_{i,p}).
%\ee
Consider the operators
\be\label{mn5}
H_i(u)=R_i(u)\frac{\prod\limits_{s\neq i}\prod\limits_{r=1}^{-a_{si}}
P_s\big(u-\ts{\i\h}{2}(\al_i+r\al_s,\al_s)\big)}
{P_i(u)P_i(u-\ts{\i\h}{2}(\al_i,\al_i))}\,,
\ee
\be\label{mn6}
E_i(u)=d_i^{-1/2}\sum\limits_{k=1}^{m_i}
\frac{\prod\limits_{s=i+1}^\ell\prod\limits_{r=1}^{-a_{si}}
P_s(\gamma_{i,k}-\ts{\i\h}{2}(\al_i+r\al_s,\al_s))}
{(u-\gamma_{i,k}) \prod\limits_{p\neq k}(\gamma_{i,k}-\gamma_{i,p})
}\,\beta_{i,k}^{-1}\,,
\ee
%where, by definition, $\prod_{s=1}^0(\ldots):=1$.
\be\label{mn7}
F_i (u)=-d_i^{-1/2}\sum\limits_{k=1}^{m_i}R_i(\gamma_{i,k}+
\ts{\imath\hbar}{2}(\alpha_i,\alpha_i))\;\times\\
\frac{\prod\limits_{s=1}^{i-1}\prod\limits_{r=1}^{-a_{si}}
P_s (\gamma_{i,k}-\ts{\i\h}{2}(\al_i+r\al_s,\al_s)+
\ts{\i\h}{2}(\al_i,\al_i))}{(u-\gamma_{i,k}-\frac{\imath\hbar}{2}
(\alpha_i,\alpha_i))\prod\limits_{p\neq k}(\gamma_{i,k}-\gamma_{i,p})}
\,\beta_{i,k}\,, \\i=1,\ldots,\ell,
\ee
where $R_i(u)$ will be specified below.

\begin{te}\label{MainTH}
\noindent (i). For any set of integer numbers $\{m_i\}$ satisfying the
condition $l_i:=\sum_{j=1}^{\ell}m_j a_{ji}\in\ZZ_+$,
 consider the polynomials
\be\label{mn4}
R_i(u)=\prod\limits_{k=1}^{l_i}(u-\nu_{i,k}),
\ee
where
$\{\nu_{i,k}\,,i\in\Gamma\,,k=1,\ldots,l_i\}$
is a set of arbitrary complex parameters.
Then the operators (\ref{mn5})-(\ref{mn7}) considered as formal power
series in $u^{-1}$, form a representation of $Y(\frak{g})$ in the space
$\cal M$. This representation is parameterized  by a choice of
$\{m_i\}$ obeying the above restrictions ($l_i\geq 0$) and by arbitrary
complex parameters $\{\nu_{i,k}\}$.

\noindent (ii). Let $\{m_i\}$ be  arbitrary integers and $R_i(u)$ be
rational functions of the following form
\be
R_i(u)=\frac{\prod\limits_{k_+=1}^{l^+_i}(u-\nu^+_{i,k_+})}
{\prod\limits_{k_-=1}^{l^-_i}(u-\nu^-_{i,k_-})},
\ee
where
$\{\nu^{\pm}_{i,k_{\pm}}\,,i\in\Gamma\,,k=1,\ldots,l^{\pm}_i\}$
is a set of arbitrary complex parameters and
$l_i^+-l_i^-=\sum_{j=1}^{\ell}m_ja_{ji}$. Then  the operators
 (\ref{mn5})-(\ref{mn6})  considered as formal power
series in $u^{-1}$, form a representation of $Y(\frak{b})$ in the space
$\cal M$. This representation is parameterized  by the choice of
$\{m_i\}$  and by arbitrary
complex parameters $\{\nu^{\pm}_{i,k}\}$.
\end{te}
{\it Proof}.
To prove the theorem introduce the following difference operators
\be\label{mn24}
\vk_{i,k}=d_i^{-1/2}\;
\frac{\prod\limits_{s=i+1}^\ell\prod\limits_{r=1}^{-a_{si}}
P_s(\gamma_{i,k}-\ts{\i\h}{2}(\al_i+r\al_s,\al_s))}
{\prod\limits_{p\neq k}(\gamma_{i,k}-\gamma_{i,p})}\,\beta_{i,k}^{-1}\,,
\ee
and
\be\label{mn25}
\hspace{-1.5cm}
\vk'_{i,k}=-d_i^{-1/2}R_i(\gamma_{ik}+\ts{\imath\hbar}{2}(\alpha_i,\alpha_i))
\frac{\prod\limits_{s=1}^{i-1}
\prod\limits_{r=1}^{-a_{si}}
P_s (\gamma_{i,k}-\ts{\i\h}{2}(\al_i+r\al_s,\al_s)+
\ts{\i\h}{2}(\al_i,\al_i))}{\prod\limits_{p\neq k}(\gamma_{i,k}-\gamma_{i,p})
}\,\beta_{i,k}\,.
\hspace{-1cm}
\ee
where $R_i(u)$ are  rational functions of $u$.
%Below $(ik)\neq (jl)$ means $i\neq j$ and $k\neq l$
It is easy to show that the operators $\gamma_{i,k},\,\vk_{i,k},\,
\vk_{i,k}'$ satisfy the relations
\be\label{mn20}
\vk_{i,k}\gamma_{j,l}-\gamma_{j,l}\vk_{i,k}=
-\ts{\imath\hbar}{2}(\alpha_i ,\alpha_i)\delta_{i,j}\delta_{k,l}\vk_{i,k}\,,\\
\vk'_{i,k}\gamma_{j,l}-\gamma_{jl}\vk'_{i,k}=\ts{\i\hbar}{2}\,
(\al_i,\al_i)\delta_{ij}\delta_{kl}\vk'_{i,k}\,,
\ee
\be\label{mn21}
(\gamma_{i,k}-\gamma_{j,l}-\ts{\imath\hbar}{2}(\alpha_i ,\alpha_j))
\vk_{i,k}\vk_{j,l}=
(\gamma_{i,k}-\gamma_{j,l}+\ts{\imath\hbar}{2}(\alpha_i ,\alpha_j))
\vk_{j,l}\vk_{i,k}\,,\\
\vk'_{i,k}\vk'_{j,l}(\gamma_{i,k}-\gamma_{j,l}+\ts{\i\hbar}{2}\,(\al_i,\al_j))
=\vk'_{j,l}\vk'_{i,k}(\gamma_{i,k}-\gamma_{j,l}-
\ts{\i\hbar}{2}\,(\al_i,\al_j))\,,
\ee
\be\label{mn15}
\hspace{1cm}
[\vk_{i,k},\vk'_{j,l}]=0,\,\,\, i\neq j,\,\,\,
[\vk_{i,k},\vk'_{i,l}]=0,\,\,\, k\neq l\,,
\ee
\be\label{mn16}
\vk_{i,k}\vk'_{i,k}=-d_i^{-1}
\frac{R_i(\gamma_{i,k})
\prod\limits_{s\neq i}\prod\limits_{r=1}^{-a_{si}}
P_s(\gamma_{i,k}-\ts{\i\h}{2}(\al_i+r\al_s,\al_s))}
{\prod\limits_{p\neq k}(\gamma_{i,k}-\gamma_{i,p})
\prod\limits_{p\neq k}(\gamma_{i,k}-\gamma_{i,p}-\ts{\i\hbar}{2}(\al_i,\al_i))}\,,\\
\hspace{-1cm}
\vk_{i,k}'\vk_{i,k}= -d_i^{-1}
\frac{R_i(\gamma_{i,k}+\ts{\i\hbar}{2}(\al_i,\al_i))
\prod\limits_{s\neq i}\prod\limits_{r=1}^{-a_{si}}
P_s(\gamma_{i,k}-\ts{\i\h}{2}(\al_i+r\al_s,\al_s)+
\ts{\i\hbar}{2}(\al_i,\al_i))}
{\prod\limits_{p\neq k}(\gamma_{i,k}-\gamma_{i,p})
\prod\limits_{p\neq k}(\gamma_{i,k}-\gamma_{i,p}+\ts{\i\hbar}{2}(\al_i,\al_i))}\,.
\hspace{-1cm}
\ee
Now let us define the generators as
\be\label{mn9}
H_i(u)=R_i(u)\,\frac{\prod\limits_{s\neq i}\prod\limits_{r=1}^{-a_{si}}
P_s\big(u-\ts{\i\h}{2}(\al_i+r\al_s,\al_s)\big)}
{P_i(u)P_i(u-\ts{\i\h}{2}(\al_i,\al_i))}\,,
\ee
\be\label{mn10}
E_i(u)=\sum\limits_{k=1}^{m_i}\frac{1}{u-\gamma_{i,k}}\,\vk_{i,k}
\,,\ \ \ i\in\Gamma,
\ee
\be\label{mn11}
F_{i}(u)=\sum\limits_{k=1}^{m_i}\vk'_{i,k}\,\frac{1}{u-\gamma_{i,k}}\,,
\ \ \ i\in\Gamma\,,
\ee
and let $R(u_i)$ be rational functions  compatible with the
expansion (\ref{int7}). Then the  relations (\ref{cr1}),(\ref{cr2}), and
(\ref{cr4}) may be derived by straightforward calculations. If we further
restrict $R_i(u)$ to be  polynomial functions then the additional relations
(\ref{cr3}) hold. To complete the proof of the theorem one should verify
the relations (\ref{cr5}), which forms in fact the only non-trivial part
of the proof.  One can see that such a verification is reduced to the
following combinatorial lemma.                                  \square

\begin{lem}
Let $\gamma_{ik}, \vk_{ik},$ and $ \vk'_{ik}$ satisfy the relations
(\ref{mn20}), (\ref{mn21}). Then for any $n=2,3,4$ the following
formulas holds:
\be\label{mn22}
\sum_{\sigma\in\frak S_n}[\vk_{i,k_{\sigma(1)}},
[\ldots,[\vk_{i,k_{\sigma(n)}},\vk_{j,l}]]...]\,=\\ =
\Big(\frac{\imath\hbar(\alpha_i ,\alpha_i)}{2}\Big)^n\;
\prod_{s=0}^{n-1}(a_{ij}+s)\cdot
\Big(\sum_{\sigma\in\frak S_n}\tilde\vk_{i,k_{\sigma(1)}}\cdot\ldots\cdot
\tilde\vk_{i,k_{\sigma(n)}}\Big)\cdot\vk_{j,l}\,,
\ee
\be\label{mn22'}
\sum_{\sigma\in\frak S_n}
[\vk'_{i,k_{\sigma(1)}},[\ldots,[\vk'_{i,k_{\sigma(n)}},
\vk'_{j,l}]]...]\,=\\ =\Big(\frac{\imath\hbar(\alpha_i ,\alpha_i)}{2}\Big)^n\;
\prod_{s=0}^{n-1}(a_{ij}+s)\cdot\Big(\sum_{\sigma\in\frak S_n}
\tilde\vk'_{i,k_{\sigma(1)}}\cdot\ldots\cdot\tilde\vk'_{i,k_{\sigma(n)}}\Big)
\cdot\vk'_{j,l}\,,
\ee
where
\be\label{mn23}
\tilde\vk_{i,k}:=
\frac{1}{\gamma_{i,k}-\gamma_{j,l}+\frac{\imath\hbar(\alpha_i ,\alpha_j)}{2}}\,
\vk_{i,k}\,,\\
\tilde\vk'_{i,k}:=\vk'_{i,k}\,
\frac{1}{\gamma_{i,k}-\gamma_{j,l}+\frac{\imath\hbar(\alpha_i,\alpha_j)}{2}}\,
\,.
\ee
\end{lem}
{\it Proof.} We outline the proof of the only non-trivial relations
(\ref{mn22}). Thus, we should calculate the following expression for
$n=2,3,4$
\bqa
X_n=\sum_{\sigma\in\frak S_n}[\vk_{i,k_{\sigma(1)}},
[\ldots,[\vk_{i,k_{\sigma(n)}},
\vk_{j,l}]]...]
\eqa
We consider only the case of non-coincident indexes $k_i\neq k_j$ for
any $i\neq j$. The proof of the general case is quite similar.
Let $\eta_{ij}=\frac{\imath\hbar(\alpha_i,\alpha_j)}{2}$.
For $n=2,3,4$,  $X_n$ may be represented as a sum of  $n!$ terms
of the first type:
\be\label{ap7}
\eta_{ij}^n\prod_{k=1}^n(\gamma_{i,k}-\gamma_{j,m})^{-1}
[\vk_{i,k_{\sigma(1)}},[\ldots,[\vk_{i,k_{\sigma(n)}},
 \vk_{j,l}]_{_+}]_{_+}...]_{_+}
\ee
and $n^n-n!$ terms of the second  type:
\be\label{ap8}
\eta_{ii}^s\eta_{ij}^{n-s}\prod_{k=1}^s(\gamma_{i,k}-\gamma_{j,m})^{-1}
\prod_{\alpha,\beta}(\gamma_{i,\alpha}-\gamma_{i,\beta})^{-1}\,\times\\

[\vk_{i,k_{\sigma(1)}},[\ldots, \vk_{j,l}
[\vk_{i,k_{\sigma(\alpha)}},\vk_{i,k_{\sigma(\beta)}}]]_{_+}]_{_+}...]_{_+}\,.
\ee One can reduce the terms of the second type to
$(\begin{smallmatrix}n\\2\end{smallmatrix})$
terms of the first type as follows. Given
$I_n:=\{1,\ldots,n\}$ and the set of variables $\{a_i,\,i\in
I\},\, d$ the following iterative formula holds.
$$\prod\limits_{i\in I_n}(a_i-d)^{-1}=\sum\limits_{r\in I_n}(a_r-d)^{-1}
\prod\limits_{i\in I_n\setminus\{r\}}(a_i-a_r)^{-1}.$$ The left
hand side $\prod\limits_{s=1}^n(a_s-d)^{-1}$ of this formula is
exactly of the first type for $a_s:=\gamma_{i,k_s}$ and
$d:=\gamma_{j,l}$. The iterations on the right hand side coincide
with all the terms of the second type (\ref{ap8}) and thus we
are left  with only the terms of the first type. The simple
transformations  then lead to  (\ref{mn22}).  \square

Finally, let us briefly explain how these representations naturally arise
from the relations (\ref{y6})-(\ref{rm13}). Due to (\ref{y6})  $A^{(s)}_i,
i\in\Gamma, s=0,1,2, \ldots$
generate a commutative subalgebra of $Y(\frak{g})$. We would like to
construct the representation in the space of functions of the finite
collection of the variables $\{\gamma_{ik}\}$ such that  $A_i^{(s)}$ act
through the multiplication by certain functions of $\{\gamma_{ik}\}$.
It is natural to look for the representation of $A_i(u)$ in the form
$A_i(u)=X_i(u)P_i(u)$ where $P_i(u)$ are
given by (\ref{mn1}) and $X_i(u)=1+\sum_{s=0}^{\infty}X_i^{(s)}
u^{-1-s}$  are  some $\gamma_{ik}$-independent series. From the
commutation relations (\ref{y6})-(\ref{rm13}) one derives that
$B_i(\gamma_{i,k})$ and $(C_i(\gamma_{i,k}))^{-1}$ are  proportional
to the shift operator (\ref{shift}). Therefore, by (\ref{rf8}) the residues
of $E_i(u)$ and $F_i(u)$ are proportional to the $B_i(\gamma_{i,k})$ and
$C_i(\gamma_{i,k})$ respectively.
This explains the ansatz (\ref{mn10}), and (\ref{mn11}) for
the generators $E_i(u)$ and $F_i(u)$.

\section{Symplectic leaves of the Yangian and the monopole moduli spaces}
In this section we describe the Poisson geometry relevant to the
description of the Yangian representations proposed above. It appears
that this leads to the direct connection with moduli spaces of
$G$-monopoles such that $Lie(G)=\frak{g}$.

Let $Y_{cl}(\frak{g})$ and $Y_{cl}(\frak{b})$ be the Poisson
algebras corresponding to the classical limit of
 $Y(\frak{g})$ and $Y(\frak{b})$ in the sense
of \cite{Dr4}, \cite{STS}. The elements $Y_{cl}(\frak{g})$ may be
described as functions on the formal loop group $LG_-$ based at
the trivial loop $g(u)=e$ where $e\in G$ is unity element.  We
use its parameterization in terms of the  infinite series of the
form $F(u)=\sum_{s=0}^{\infty}F^{(s)}u^{-s-1}$.

The description of  the Poisson algebra
$Y_{cl}(\frak{g})$ in terms of the generators and relations
can be obtained from (\ref{cr1})-(\ref{cr5}) by taking the limit $\hbar
\rightarrow 0$:
\be\label{cr1'}
\{h_i(u),h_j(v)\}=0\;,
\ee
\be
\{h_i(u),e_j(v)\}=-\, (\alpha_i,\alpha_j)\,
\frac{h_i(u)(e_j(u)-e_j(v))}{u-v}\,,\\
\{h_i(u),f_j(v)\}= (\alpha_i,\alpha_j)\,
\frac{h_i(u)(f_j(u)-f_j(v))}{u-v}\,,
\ee
\be\label{cr3'}
\{e_i(u),f_j(v)\}=-\,\delta_{ij}\frac{h_i(u)-h_i(v)}
{u-v}\;,
\ee

\be\label{cr4'}
\hspace{-0.5cm}
\{e_i(u),e_i(v)\}=
-(\alpha_i,\alpha_i)\frac{(e_i(u)-e_i(v))^2}{u-v},\\
\hspace{-0.5cm}
\{f_i(u),f_i(v)\}=
(\alpha_i,\alpha_i)\frac{(e_i(u)-e_i(v))^2}{u-v},\\
\hspace{-0.5cm}
\{e_i(u),e_j(v)\}=-
(\alpha_i,\alpha_j)\frac{e_i(u)(e_j(u)-e_j(v))}{u-v}-
\frac{\{e_i^{(0)},(e_j(u)-e_j(v))\}}{u-v},\\
\hspace{-0.5cm}
\{f_i(u),f_j(v)\}=
(\alpha_i,\alpha_j)\frac{f_i(u)(f_j(u)-f_j(v))}{u-v}-
\frac{\{f_i^{(0)},(f_j (u)-f_j(v))\}}{u-v},\\
i\neq j,\,\,\, a_{ij}\neq 0;
\ee
\be\label{cr5'}
\sum_{\sigma\in\frak S_n} \{e_i(u_{\sigma(1)}),\{e_i(u_{\sigma(2)}),
\ldots,\{e_i(u_{\sigma(n)}),e_j(v)\}\ldots\}\}=0\,,\\
\sum_{\sigma\in\frak S_n}\{f_i(u_{\sigma(1)}),\{f_i(u_{\sigma(2)}),
\ldots,\{f_i(u_{\sigma(n)}),f_j(v)\}\ldots\}\}=0\;,
\\ n=1-a_{ij},\ \ (i\neq j)\,.
\ee
There exists the following simple  interpretation of the generators
$h_i(u),e_i(u),f_i(u)$. Let us fix the Gauss
decomposition of an element $g(u)\in LG_-$:
\bqa \label{Gauss}
g(u)=\exp\big(\sum_{\a}f_{\a}(u)F_{\a}\big)\cdot
\exp\big(\sum_{i=1}^\ell \phi_i(u)H_i\big)\cdot
\exp\big(\sum_{\a}e_{\a}(u)E_{\a}\big)\,,
\eqa
where $H_i,F_{\a},E_{\a}$ provide  a basis of $\mathfrak{g}$
labeled by  positive roots $\a$, and $\phi_i(u),f_{\a}(u),e_{\a}(u)$ are
the local exponential coordinates on
the group $LG_-$.  Note that we consider the functions on the formal
loop group  $LG_-$ (i.e.  we deal with functions on the formal
neighbourhood of $e\in LG_{-}$) and thus  the Gauss decomposition (\ref{Gauss})
is valid "everywhere".
The  functions $e_i(u):=e_{\a_i}(u)$, $f_i(u):=f_{\a_i}(u)$
corresponding to the simple roots $\a_i$  together with
$h_i(u)=\exp\{-\sum_{j=1}^\ell a_{ji}\phi_j(u)\}$ give us a
set of the generators satisfying the relations (\ref{cr1'})-(\ref{cr5'}).
The local coordinates $\phi_i(u),e_{i}(u),f_i(u)$
may be  expressed  explicitly  in terms of the matrix elements of the
fundamental representations of $U(\frak{g})$.
Let $\{\pi_i\}$ be a set of  fundamental  representations corresponding
to the fundamental weights $\{\omega_i\}$ of
$\mathfrak{g}$ and $v^{(i)}_{+/-}$ be the highest/lowest vectors in
these representation. Denote by
$a_i(u)$, $b_i(u)$, $c_i(u)$, $d_i(u)$ the following formal series
\be\label{matrixel}
a_i(u)=\<v_-^{(i)}|\pi_{i}(g(u))|v_+^{(i)}\>\,,\\
b_i(u)=\<v_-^{(i)}|\pi_{i}(g(u))\pi_{i}(F_i)|v_+^{(i)}\>\,,\\
c_i(u)=\<v_-^{(i)}|\pi_{i}(E_i)\pi_{i}(g(u))|v_+^{(i)}\>\,,\\
d_i(u)=\<v_-^{(i)}|\pi_{i}(E_i)\pi_{i}(g(u))\pi_{i}(F_i)|v_+^{(i)}\>\,.
\ee
\begin{lem}\label{matrixelone}
The  coordinates corresponding to the simple roots and Cartan elements
 entering the Gauss decomposition (\ref{Gauss})
can be expressed through the matrix elements (\ref{matrixel}) as follows
\be
e^{\phi_i(u)}=a_i(u)\,,\\
e_{i}(u)=\frac{b_i(u)}{a_i(u)}\,,\\
f_{i}(u)=\frac{c_i(u)}{a_i(u)}\,.
\ee
\end{lem}
The variables $a_i(u),b_i(u),c_i(u),d_i(u)$ are the classical
counterparts of the variables $A_i(u)$, $B_i(u)$, $C_i(u)$, $D_i(u)$
defined in (\ref{rf3}),(\ref{rf8}),(\ref{rf12}).
A similar description of $Y(\frak{gl}(\ell+1))$ was given in \cite{Dr3}
(see also \cite{FJ}, \cite{Io}). In this case the functions
$a_i(u),b_i(u),c_i(u),d_i(u)$ are given by the minors of the matrix
$\pi_{*}(g(u))$ where $\pi_*$ is the tautological representation
$\pi_*:\frak{gl}(\ell+1)\rightarrow{\rm End}(\CC^{\ell+1})$, which
agrees with (\ref{matrixel}).  Let us remark that
  (\ref{matrixel}) may be considered as a classical counterpart of the universal R-matrix map of
the  dual Hopf algebra with the opposite comultiplication $A^0$
to the Hopf algebra $A$  (see
\cite{Dr2} for details).
Together with the explicit Gauss form   of the universal $R$-matrix
this should provide  the quantum version of Lemma \ref{matrixelone}.
The simplest example of  $Y(\mathfrak{sl}(2))$ may be extracted from
\cite{KT}.

Now consider the classical counterpart of the Yangian representations
constructed in the Section 3. These representations have the following
property: the images of $E_i(u)$   are rational operator-valued
functions in $u$ with simple poles. Given the Poisson brackets
 (\ref{cr1'})-(\ref{cr5'}) on
$LG_-$, the representations of the Yangian correspond to the symplectic
leaves in $LG_-$. A similar description holds for $Y_{cl}(\frak{b})$
 in terms of the symplectic leaves in $LB_-$.
We define the
symplectic leaf  $\mathcal{O}$ to be
rational if the restriction of the generators $e_i(u)$
 is  a rational function over $u$ and let
 $\mathcal{O}^{(0)}\subset \mathcal{O}$ be an open part corresponding
to $e_i(u)$ having only simple poles.  Thus the symplectic leaves
corresponding to  the representations constructed in Section 3 are
rational. One can describe the symplectic leaves as symplectic
manifolds as follows.
 Open parts $\mathcal{O}^{(0)}$ of the rational
  symplectic leaves  in  $LB_-$
corresponding to the representations constructed in Theorem
\ref{MainTH} are isomorphic (as abstract manifolds) to the open
subsets of the space of $\widetilde{\mathcal{M}}_b({\bf m })$ of the
based rational maps \bqa\label{maps} e=(e_1,\cdots
,e_{\ell}):(\mathbb{P}^1,\infty)
\stackrel{e}{\longrightarrow}(\underbrace{\mathbb{P}^1\times
\cdots \times \mathbb{P}^1}_\ell\,,0 \times \cdots \times 0), \eqa
of the fixed multi-degree ${\bf m }=(m_1,\cdots ,m_{\ell})$  where $e_i(u)$
are the generators of  $Y_{cl}(\frak{b})$ corresponding to simple
roots. Analogously open parts $\mathcal{O}^{(0)}$ of the rational
symplectic leaves of $LG_-$ corresponding to the representations
constructed in Theorem \ref{MainTH} are isomorphic to the open
subsets of $\widetilde{\mathcal{M}}_b({\bf m})$  with additional
restrictions $\sum_{j=1}^{\ell}m_ja_{ji}=l_i \in \ZZ_+.$

Taking into account the results of \cite{FM} one can reformulate
the description of the symplectic leaves as follows.
   Consider the space $\mathcal{M}({ \bf m})$
of the holomorphic maps $\mathbb{P}^1\rightarrow G/B$ of
multi-degree ${\bf m}=(m_1,\ldots,m_\ell)\in \Lambda_W^{\vee}$, where
$\Lambda_W^{\vee}=H_2(G/B,\mathbb{Z})$ is the co-weight lattice of
$\frak{g}$. It will be useful to consider $G/B$ as a manifold
parameterizing the Borel subgroups in $G$. Choose  some
Borel subgroup $B_+$  and let $b_+\subset G/B$ be the corresponding
point in the flag manifold.  Let us fix the local coordinate
 on $\mathbb{P}^1$ and consider the evaluation map
$ev_\infty:\,\mathcal{M}({\bf m})\longrightarrow G/B$ defined as $ev_{
\infty}:\, f \rightarrow f(\infty)$.  Thus $\mathcal{M}({\bf m})$
 is supplied with the structure of the fibred space
over $G/B$ and the fibre  is naturally identified with the
 moduli space  $\mathcal{M}_b({\bf m})$ of the based holomorphic
 maps $f:(\mathbb{P}^1,\infty) \rightarrow (G/B,b_+)$  of the multi-degree
${\bf m}$. It appears that the open part of
 $\mathcal{M}_b({\bf m})$  can be naturally identified with the moduli space
$\widetilde{\mathcal{M}}_b({\bf m})$ introduced above. Actually this follows  from the results
of  Drinfeld in the form presented  in \cite{FM}. Thus it was shown in
\cite{FM} that $\mathcal{M}_b({\bf m})$ is  a smooth manifold
of dimension
$${\rm dim}\,\mathcal{M}_b({\bf m})=2|{\bf m}|=2(m_1+\ldots+m_\ell).$$
The explicit description of the manifold  $\mathcal{M}_b({\bf m})$
can be obtained by generalizing the classical Pl\"{u}cker  embedding of $G/B$ into the
product  $\prod_{i \subset
\Gamma }\mathbb{P}(V_{\omega_i})$  of  the projectivisations of
the fundamental representations $V_{\omega_i}$ as follows.
  Let $\pi_{\l}:U(\mathfrak{g})\rightarrow End(V_{\l})$ be the
  irreducible  representation of the universal enveloping algebra
$U(\mathfrak{g})$ with  the highest weights $\l$ and
$\mathcal{V}_{\l}=V_{\l}\otimes \mathcal{O}_{\mathbb{P}^1}$ be the
corresponding trivial vector bundles on $\mathbb{P}^1$.
Using the local coordinate we identify
$\Gamma(\mathbb{A}^1,\mathcal{V}_{\l})=V_{\l}\otimes \mathbb{C}[u]$.
Denote by $v_+^{\l}$ the highest weight vectors in $V_{\l}$  with respect to the
Borel subgroup $B_+$. Similarly let $v_-^{\l}$ be  the
lowest weight vector in the dual representation $V^{\vee}_{\l}$
normalized by the condition  $\<v_-^{\l}|v_+^{\l}\>=1$. We also introduce the additional
set of vectors $\pi_{\l}(F_{i})v_+^{\l}$, $\pi_{\l}( E_i)v_-^{\l}$
where $F_{i}$ and $E_{i}$  are the generators corresponding  to
the simple roots $\a_i$. Thus given a section
$v^{\l} \in
\Gamma(\mathbb{A}^1,\mathcal{V}_{\l})$
we have the following decomposition $v^{\l}(u)=a_{\l}(u)\cdot
v_+^{\l}+\sum_{i\in \Gamma}b^i_{\l}(u) \pi_{\l}(F_i)v_+^{\l}+\phi^{\l}(u)$ with
$\phi^{\l}\subset \Gamma(\mathbb{A}^1,\mathcal{V}_{\l})$
satisfying $\<v_-^{\l}|\phi^{\l}\>=
\<v_-^{\l}|\pi_{\l}(E_i)|\phi^{\l}\>=0$.

According to Drinfeld (see \cite{FM} for the details)
 the moduli space $\mathcal{M}_b({\bf m})$ is isomorphic to the space $Z_{\bf m}$
of the collections of sections $v^{\l}(u)$
for each $\l \in \Lambda^+_W$ satisfying the conditions
\begin{enumerate}
\item The polynomial $a_{\l}(u)$ is monic of degree $\< m,\l\>$;
\item The degree of $(v^{\l}-a_{\l}(u)v_+^{\l})$ is strictly less then $\< m,\l \>$;
\item For any $G$-equivariant morphism $\phi: V_{\l}\otimes V_{\mu}\rightarrow
  V_{\nu}$ such that $\nu=\mu+\lambda$ and the conjugated morphism
  satisfies   $\phi^*(v_-^{\nu})=v_-^{\l}\otimes v_-^{\mu}$  we have
  $\phi(v^{\l}\otimes v^{\mu})=v^{\nu}$;
\item For any $G$-equivariant morphism $\phi: V_{\l}\otimes V_{\mu}\rightarrow
  V_{\nu}$ such that $\nu<\mu+\lambda$   we have
  $\phi(v^{\l}\otimes v^{\mu})=0$.
\end{enumerate}
It is easy to see that the set  $\{v^{\l}, \lambda\in\Lambda^+_W\}$
 satisfying these conditions is
determined by its subset  $\{v^{\omega_i}\}$ corresponding to
fundamental representations $\omega_i$.  Moreover,
given  arbitrary polynomials $a_{\omega_i}(u)$ and $b^i_{\omega_i}(u)$
satisfying the conditions (1) and (2) above (i.e. $a_{\omega_i}(u)$ are
monic and $deg(a_{\omega_i})=deg(b^i_{\omega_i})+1=m_i$)
 there exist such $\phi^{\omega_i}(u)$ that for
$v^{\omega_i}(u)=a_{\omega_i}(u)\cdot
v_+^{\omega_i}+b^i_{\omega_i}(u) \pi_{\omega_i}(F_i)v_+^{\omega_i}+
\phi^{\omega_i}(u)$
the conditions (3) and (4) hold (see \cite{FM}, \cite{FKMM} for
details).  Let us consider the subset of the polynomials $a_{\omega_i}$
and $b_{\omega_i}$  such that the roots $\gamma_{i,k}$ of  $a_{\omega_i}(u)$ do not
coincide $\gamma_{i,k}\neq \gamma_{j,l}$ for $(i,k)\neq (j,l)$.
The space  of such  polynomials  $a_i(u)\equiv a_{\omega_i}(u)$,
$b_i(u)\equiv b^i_{\omega_i}(u)$
is $2|{\bf m}|$-dimensional and thus is
isomorphic to the open subspace in the moduli space
$\mathcal{M}_b({\bf m })$. Note that the polynomials $a_i(u)$ and $b_i(u)$   define the   map
$e\in\widetilde{\mathcal{M}}_b({\bf m})$
$$(\mathbb{P}^1,\infty)\stackrel{e}{\longrightarrow}(\underbrace
{\mathbb{P}^1\times\ldots\times\mathbb{P}^1}_\ell\,,0\times\ldots\times 0),$$
  given by
$$e(u)=(b_1(u)/a_1(u))\times\ldots\times(b_\ell(u)/a_\ell(u)).$$
Therefore we have established  the isomorphism of the open parts of the
moduli spaces
$$\phi:\mathcal{M}_b({\bf
  m})\longrightarrow\widetilde{\mathcal{M}}_b({\bf m}).$$
One can summarize this  in the following
\begin{prop}
{\it (i) The open parts $\mathcal{O}^{(0)}$ of the rational
symplectic leaves  of $Y_{cl}(\frak{b})$ corresponding to the
representations constructed in
  Theorem
\ref{MainTH} are
 isomorphic to  the  open parts  of the spaces of the based  maps
\bqa \label{maps1}
 (\mathbb{P}^1,\infty)\rightarrow (G/B,b_+)\eqa
of the fixed multi-degree ${\bf m}=(m_1,\cdots ,m_{\ell})\in H_2(G/B,\ZZ)$.

(ii) The open parts $\mathcal{O}^{(0)}$ of the rational symplectic
leaves of $Y_{cl}(\frak{g})$ corresponding to the representations
constructed in Theorem \ref{MainTH} are isomorphic to the  spaces
of the based  maps (\ref{maps1}) with additional restrictions
$\sum_{j=1}^{\ell}m_ja_{ji}=l_i \in \ZZ_+. $}
\end{prop}
The connection with the explicit parameterization used in the previous
sections is as follows.   Let us  parameterize the open subset
$U\subset\widetilde{\mathcal{M}}_b({\bf m})$ by the following \'etale
coordinates
$$(x_{i,k},\,y_{i,k}),\quad
i=1,\ldots,\ell,\quad k=1,\ldots,m_i
$$
defined by the conditions
$$a_i(x_{i,k})=0,\quad
y_{i,k}=b_i(x_{i,k}). $$ Then the  coordinates $(x_{ik},y_{ik})$  are
related to the coordinates $(\gamma_{i,k},\vk_{i,k})$  by the simple redefinition
\bqa
x_{i,k}&=&\gamma_{i,k}\,,\\
y_{i,k}&=&\vk_{i,k}\prod_{s\neq k}(\gamma_{i,k}-\gamma_{i,s})\,.
\eqa
Note that the classical limit of the relations
(\ref{mn20}),(\ref{mn21}) provides the open part of the space
$\mathcal{M}_b(m)$ with the holomorphic symplectic structure
\bqa
\{ \gamma_{i,k},\gamma_{j,l}\}&=&0,\\
\{ \gamma_{i,k},\vk_{j,l}\}&=&
\ts{1}{2}(\alpha_i,\alpha_i)\delta_{i,j}\delta_{k,l}\vk_{j,l}\,,\\
\{\vk_{i,k},\vk_{j,l}\} &=& (\alpha_i ,\alpha_j)
  \frac{\vk_{i,k}\vk_{j,l}}{\gamma_{i,k}-\gamma_{j,l}}\,, \,\,\,\ (i,k)\neq
(j,l), \eqa or equivalently in the coordinates $(x_{i,k},y_{i,k})$
\bqa \label{ssF}
 \{ x_{i,k},x_{j,l}\}&=&0,\\
 \{y_{i,k},y_{i,l} \}&=&0,\\
 \{x_{i,k},y_{j,l}\}&=&
\ts{1}{2}(\alpha_i,\alpha_i)\delta_{i,j}\delta_{k,l}
y_{j,l},\\
\{y_{i,k},y_{j,l}\} &=& (\alpha_i ,\alpha_j)
  \frac{y_{i,k}y_{j,l}}{x_{i,k}-x_{j,l}}\,, \,\,\,\  i\neq j.
\eqa
 In general the symplectic leaves of the Poisson-Lie
group $G$ are connected components of the intersections of the
double cosets of the Poisson-Lie dual $G^*$ in $G\times G$ with
the diagonal $G\subset G\times G$ \cite{STS}. We are going to
discuss the connection of this description with the
algebro-geometric description considered above in the separate
publication.

It appears that description of the symplectic leaves  of the
Poisson-Lie groups  associated with the Yangian   given  in this
  section provides a direct connection
with the moduli spaces of the $G$-monopoles with the maximal symmetry
breaking. These moduli spaces  are also given by the spaces  of
the based maps $\mathbb{P}^1\rightarrow G/B$
\cite{Hurtubise}, \cite{HM}, \cite{Ja}.
 Our construction provides the holomorphic symplectic
structure on these spaces.
The explicit  description of the  holomorphic symplectic structure
on the moduli space of the monopoles was given in the case of
$G=SU(N)$ in \cite{B2} (generalizing the results for $SU(2)$ of \cite{AH})
and for general case in \cite{FKMM}. It turns out  that our description
of the  coordinates on the moduli space and the expression for the
Poisson structure in coordinates $(x_{i,k},y_{i,k})$
exactly matches the
description given in \cite{FKMM}.
  Thus  the representation
constructed in Section 2 can be considered as a quantization of the moduli
space of the monopoles.

\end{document}